\documentclass[a4,reqno,11pt]{amsart}
\usepackage[margin=30mm]{geometry}
\usepackage{amssymb,amsmath,color}
\usepackage{hyperref}
\newcommand{\qbin}[2]{\genfrac{[}{]}{0pt}{}{#1}{#2}_q}

\theoremstyle{definition}

\begin{document}

\title{A $q$-analogue of the Koecher-Leshchiner\\
generating function of odd zeta values}

\author[R. Tauraso]{Roberto Tauraso}

\address{Dipartimento di Matematica, % \\
Università di Roma ``Tor Vergata'', 00133 Roma, Italy}
\email{tauraso@mat.uniroma2.it}

\subjclass{Primary 11B65, 05A19; Secondary 05A15, 05A19, 05A30, 11M06.}

\keywords{$q$-analogue, zeta values, Ap\'ery-like series, central binomial coefficients, generating functions}

\begin{abstract} In the 1980s, Koecher and, independently, Leshchiner found an elegant formula for the generating function of odd zeta values. In this short note, we derive a $q$-analogue of this formula, which provides a $q$-version of the accelerated series for $\zeta(3)$ used by Ap\'ery in his famous proof of irrationality.
\end{abstract}

\maketitle

\section{Introduction}

In June 1978, during a talk, the French mathematician Roger Ap\'ery   outlined a rather unexpected proof establishing that $\zeta(3)$ is an irrational number  \cite{Apery79} (see also the survey paper \cite{Poorten79}). 
A fundamental component of his approach was identity 
\begin{equation}\label{MarkovApery}
\zeta(3)=\frac{5}{2}\sum_{k=1}^\infty\frac{(-1)^{k-1}}{k^3\binom{2k}{k}}.
\end{equation}
This accelerated series for $\zeta(3)$ was actually first found by Markov in 1890. It follows from \cite[(14)]{Markov1890}, which is rewritten below in a compact mode, by letting $a=0$,
\begin{equation*}%\label{MarkovG}
\sum_{k=a+1}^{\infty}\frac{1}{k^3}=
\sum_{k=1}^\infty\frac{(-1)^{k-1}(\frac{5}{2}k^2+3ak+a^2)}{k^5\binom{2k}{k}\binom{a+k}{a}^4}
\end{equation*}
After Ap\'ery's proof, several other related identities were either established or revived. Here, we will focus on a generating function identity for odd zeta values obtained by
Koecher \cite{Koecher80} and, independently, by Leshchiner \cite{Leshchiner81}
(see also \cite[(1-2)]{Rivoal04}): 
\begin{equation}\label{KL}
\sum_{r=0}^{\infty}
\zeta(2r+3)x^{2r}=\sum_{k=1}^{\infty}\frac{1}{k(k^2-x^2)}=
\frac{1}{2}\sum_{k=1}^{\infty}\frac{(-1)^{k-1}}{k^3\binom{2k}{k}}
\cdot \frac{5k^2-x^2}{k^2-x^2}
\cdot\prod_{m=1}^{k-1}\Big(1-\frac{x^2}{m^2}\Big),
\end{equation}
which implies, for $r\geq 0$,
$$\zeta(2r+3)=
\sum_{k=1}^{\infty} \frac{(-1)^{k-1-r}}{k^3\binom{2k}{k}}
\Big(\frac{5}{2}H_{k-1}(\{2\}^{r})+2
\sum_{j=1}^{r} \frac{(-1)^j}{k^{2j}}H_{k-1}(\{2\}^{r-j})\Big)$$
where 
$H_{k}(\{2\}^s)=\sum_{1\leq j_1<\dots<j_s\leq k}1/(j_1^2\cdots j_s^2)$.

\noindent For even zeta values, Bailey, Borwein, and Bradley using an
experimental approach discovered and then formally proved \cite[Theorem 1.1]{BBB06} that
$$\sum_{r=0}^{\infty}
\zeta(2r+2)x^{2r}=\sum_{k=1}^{\infty}
\frac{1}{k^2-x^2}
=\sum_{k=1}^{\infty}
\frac{3}{k^2\binom{2k}{k}}
\cdot \frac{\prod_{m=1}^{k-1}\big(1-\frac{4x^2}{m^2}\big)}{\prod_{m=1}^{k}\big(1-\frac{x^2}{m^2}\big)}.$$
Later in \cite[Theorem 1]{Hessami11}), Kh. Hessami Pilehrood and T. Hessami Pilehrood proved the interesting $q$-analogue, 
\begin{equation}\label{qeven}
\sum_{k=1}^{\infty}
\frac{q^k}{[k]^2_q-x^2q^{2k}}
=\sum_{k=1}^{\infty}
\frac{q^{k^2}(1+2q^k)}{[k]^2_q\qbin{2k}{k}}
\cdot \frac{\prod_{m=1}^{k-1}\big(1-\frac{x^2(1+q^m)^2}{[m]_q^2}\big)}{\prod_{m=1}^{k}\big(1-\frac{x^2q^{2m}}{[m]_q^2}\big)}.
\end{equation}
In particular, for $x=0$, we get
\begin{equation*}%\label{zeta2}
\sum_{k=1}^{\infty}\frac{q^{k^2}(1+2q^k)}{[k]^2_q\qbin{2k}{k}}
=\sum_{k=1}^{\infty}\frac{q^k}{[k]^2_q}.
\end{equation*}
Is there a $q$-analogue for the generating fuction of the odd zeta values \eqref{KL}? \medskip

\noindent In \cite[(19) and (20)]{Hessami11} the authors gave some $q$-version for the Markov-Ap\'ery formula \eqref{MarkovApery} for $\zeta(3)$. In this short note we provide a general $q$-version of \eqref{KL}:
\begin{align*}%\label{qKL}
\sum_{k=1}^{\infty}\frac{q^k+q^{2k}}{[k]_q([k]^2_q-x^2q^k)}=
\sum_{k=1}^{\infty}&\frac{(-1)^{k-1}q^{\binom{k+1}{2}}}{[k]^3_q\qbin{2k}{k}}\\
&\cdot \frac{(1+3q^k+q^{2k})[k]^2_q-x^2q^{2k}}{[k]^2_q-x^2q^k}\cdot\prod_{m=1}^{k-1}\left(1-\frac{x^2q^m}{[m]^2_q}\right).
\end{align*}
which is equivalent, after expanding the series with respect to the powers $x^{2r}$, to the family of identities
\begin{align}\label{qMarkovAperyr}
\sum_{k=1}^{\infty}\frac{q^{(r+1)k}+q^{(r+2)k}}{[k]^{2r+3}_q}
&=\sum_{k=1}^{\infty}\frac{(-1)^{k-1+r}q^{{\binom{k+1}{2}}}}{[k]^3_q\qbin{2k}{k}}\\
\nonumber&\cdot\Big((1+3q^k+q^{2k})H_{q,k-1}(\{2\}^r)
+(1+q^k)^2\sum_{j=1}^r\frac{(-1)^jq^{kj}}{[k]^{2j}_q}
H_{q,k-1}(\{2\}^{r-j})
\Big)
\end{align}
where for $s\geq 1$,  $H_{q,k}(\{2\}^s)$ is the multiple harmonic $q$-sum
$$\sum_{1\leq j_1<\dots<j_s\leq k}\frac{q^{j_1+\dots+j_s}}{[j_1]_q^2\cdots[j_s]_q^2}$$
and $H_{q,k}(\{2\}^0)=1$. It is worth noting that for $r=0$, the $q$-series \eqref{qMarkovAperyr} yields a $q$-analogue of \eqref{MarkovApery},
\begin{align*}%\label{qMarkovApery}
\sum_{k=1}^{\infty}\frac{q^k+q^{2k}}{[k]^3_q}&=\sum_{k=1}^{\infty}\frac{(-1)^{k-1}q^{{\binom{k+1}{2}}}(1+3q^k+q^{2k})}{[k]^3_q\qbin{2k}{k}}.
\end{align*}

\section{Proof of \eqref{qMarkovAperyr}}
We first introduce some notations. The $q$-analogue of a positive integer $n$ is
$$[n]_q=\frac{1-q^n}{1-q}=\sum_{j=0}^{n-1}q^j,$$
and the $q$-factorial $[n]_q!$ is defined as $\prod_{j=1}^{n}[j]_q$.
The $q$-binomial coefficient is given by 
$$\qbin{n}{k}=\frac{[n]_q!}{[k]_q![n-k]_q!}.$$
We will show the following finite form of \eqref{qMarkovAperyr}:
\begin{align*}
\sum_{k=1}^{N}\frac{q^{(r+1)k}+q^{(r+2)k}}{[k]^{2r+3}_q}
&=\sum_{k=1}^{N}\frac{(-1)^{k-1+r}q^{{\binom{k+1}{2}}}}{[k]^3_q\qbin{2k}{k}}\\
&\cdot\Big((1+3q^k+q^{2k})H_{q,k-1}(\{2\}^{r})
+(1+q^k)^2\sum_{j=1}^{r}\frac{(-1)^jq^{kj}}{[k]^{2j}_q}
H_{q,k-1}(\{2\}^{r-j})
\Big)\\
&-\sum_{k=1}^{N}\frac{(-1)^{k-1+r}q^{-\binom{k}{2}+Nk+k}}{[k]_q^3\qbin{N}{k}\qbin{N+k}{k}}H_{q,k-1}(\{2^{r}\}).
\end{align*}
Then the claim is easily obtained by letting $N\to \infty$.

\noindent The definition of $H_{q,k}(\{2^s\})$ leads to the recurrence relation
$$H_{q,k}(\{2^s\})-H_{q,k-1}(\{2^s\})=\frac{q^k}{[k]_q^2}H_{q,k-1}(\{2^{s-1}\})$$
which implies
$$q^{n-k}H_{q,k}(\{2^s\})+\frac{[n+k]_q[n-k]_q }{[k]_q^2}H_{q,k-1}(\{2^s\})=\frac{[n]_q^2}{[k]_q^2}H_{q,k-1}(\{2^s\})+\frac{q^n }{[k]_q^2}H_{q,k-1}(\{2^{s-1}\}).$$
We multiply both sides of the above equality by
$$(-1)^kq^{-\binom{n-k+1}{2}}\frac{([k]_q!)^2}{\prod_{j=n-k}^{n+k}[j]_q}.$$
Then the left-hand side becomes
$$\text{LHS}=(-1)^kq^{-\binom{n-k}{2}}\frac{([k]_q!)^2}{\prod_{j=n-k}^{n+k}[j]_q}H_{q,k}(\{2^s\})
-(-1)^{k-1}q^{-\binom{n-(k-1)}{2}}\frac{([k-1]_q!)^2}{\prod_{j=n-(k-1)}^{n+k-1}[j]_q}
H_{q,k-1}(\{2^s\})$$
while the right-hand side equals
$$\text{RHS}=(-1)^kq^{-\binom{n-k+1}{2}}\frac{([k-1]_q!)^2}{\prod_{j=n-k}^{n+k}[j]_q}
\Big([n]_q^2 H_{q,k-1}(\{2^s\})+q^n H_{q,k-1}(\{2^{s-1}\})\Big).$$
Summing both sides for $k=1,\dots,n-1$, the LHS telescopes and simplifies to
$$\text{LHS}
=\frac{(-1)^{n-1}(1+q^n)}{[n]_q\qbin{2n}{n}}
H_{q,n-1}(\{2^s\})
-\frac{q^{-\binom{n}{2}}}{[n]}H_{q,0}(\{2^s\})
$$
(note that $[2n]_q=(1+q^n)[n]_q$). Conversely, on the other side, we find
$$\text{RHS}
=\sum_{k=1}^{n-1}(-1)^kq^{-\binom{n-k+1}{2}}\frac{([k-1]_q!)^2}{\prod_{j=n-k}^{n+k}[j]_q}
\Big([n]_q^2 H_{q,k-1}(\{2^s\})+q^n H_{q,k-1}(\{2^{s-1}\})\Big).
$$
After multiplying both sides by 
$$\frac{(-1)^sq^{\binom{n}{2}+(r-s)n}(1+q^n)}{[n]_q^{2r-2s}},$$ 
we get
$$\text{LHS}
=\frac{(-1)^{s+n-1}q^{\binom{n}{2}+(r-s)n}(1+q^n)^2}{[n]_q^{2r-2s+1}\qbin{2n}{n}}
H_{q,n-1}(\{2^s\})
-\frac{q^{(r-s)n}+q^{(r-s+1)n}}{[n]^{2r-2s+1}}H_{q,0}(\{2^s\})
$$
and
\begin{align*}\text{RHS}
=q^{\binom{n}{2}}&(1+q^n)
\sum_{k=1}^{n-1}(-1)^kq^{-\binom{n-k+1}{2}}\frac{([k-1]_q!)^2}{\prod_{j=n-k}^{n+k}[j]_q}\\
&
\cdot\left( \frac{(-1)^sq^{(r-s)n}}{[n]_q^{2r-2s-2}}H_{q,k-1}(\{2^s\})-\frac{(-1)^{s-1}q^{(r-s+1)n}}{[n]_q^{2r-2s}} H_{q,k-1}(\{2^{s-1}\})\right).
\end{align*}
The next step is to sum both sides for $s=0,\dots,r-1$. We obtain
\begin{align*}\text{LHS}
&=\sum_{s=0}^{r-1}\frac{(-1)^{n-1+s}q^{\binom{n}{2}+(r-s)n}(1+q^n)^2}{[n]_q^{2r-2s+1}\qbin{2n}{n}}
H_{q,n-1}(\{2^s\})
-\frac{q^{rn}+q^{(r+1)n}}{[n]_q^{2r+1}}\\
&=\frac{(-1)^{r+n-1}q^{\binom{n}{2}}(1+q^n)^2}{\qbin{2n}{n}}
\sum_{j=1}^{r}
\frac{(-1)^jq^{nj}}{[n]_q^{2j+1}}H_{q,n-1}(\{2^{r-j}\})
-\frac{q^{rn}+q^{(r+1)n}}{[n]_q^{2r+1}},
\end{align*}
and the RHS telescopes,
\begin{align*}\text{RHS}
=q^{\binom{n}{2}}&(1+q^n)
\sum_{k=1}^{n-1}(-1)^kq^{-\binom{n-k+1}{2}}\frac{([k-1]_q!)^2}{\prod_{j=n-k}^{n+k}[j]_q}\\
&
\cdot\sum_{s=0}^{r-1}\left( \frac{(-1)^s q^{(r-s)n}}{[n]_q^{2r-2s-2}}H_{q,k-1}(\{2^s\})-\frac{(-1)^{s-1}q^{(r-s+1)n}}{[n]_q^{2r-2s}} H_{q,k-1}(\{2^{s-1}\})\right)\\
&=(-1)^{r-1}q^{\binom{n}{2}}(1+q^n)
\sum_{k=1}^{n-1}(-1)^kq^{-\binom{n-k+1}{2}}\frac{([k-1]_q!)^2}{\prod_{j=n-k}^{n+k}[j]_q}
q^{n}H_{q,k-1}(\{2^{r-1}\})\\
&=(-1)^{r-1}(1+q^n)
\sum_{k=1}^{n-1}(-1)^kq^{-\binom{k}{2}+nk}\frac{([k-1]_q!)^2}{\prod_{j=n-k}^{n+k}[j]_q}
 H_{q,k-1}(\{2^{r-1}\}).
\end{align*}
Then we proceed by summing both sides for $n=1,\dots,N$. 
After interchanging the order of summation, the RHS can be written as 
\begin{align*}
\text{RHS}&=(-1)^{r-1}\sum_{n=1}^{N}(1+q^{n})\sum_{k=1}^{n-1}
(-1)^kq^{-\binom{k}{2}+nk}\frac{([k-1]_q!)^2}{\prod_{j=n-k}^{n+k}[j]_q}H_{q,k-1}(\{2^{r-1}\})\\
&=(-1)^{r-1}\sum_{k=1}^{N-1}\frac{(-1)^kq^{-\binom{k}{2}}}{[k]_q^2\qbin{2k}{k}}
H_{q,k-1}(\{2^{r-1}\})\cdot [2k]_q!\sum_{n=k+1}^{N}
\frac{(1+q^{n})q^{nk}}{\prod_{j=n-k}^{n+k}[j]_q}\\
&=(-1)^{r-1}\sum_{k=1}^{N-1}\frac{(-1)^kq^{-\binom{k}{2}}}{[k]_q^2\qbin{2k}{k}}
H_{q,k-1}(\{2^{r-1}\})\cdot [2k]_q!
\sum_{n=k+1}^{N}\frac{1}{[k]_q}\left(
\frac{q^{nk}}{\prod_{j=n-k}^{n+k-1}[j]_q}
-\frac{q^{(n+1)k}}{\prod_{j=n+1-k}^{n+1+k-1}[j]_q}
\right),
\end{align*}
so that the RHS telescopes again
\begin{align*}
\text{RHS}&=(-1)^{r-1}\sum_{k=1}^{N-1}\frac{(-1)^kq^{-\binom{k}{2}}}{[k]_q^3\qbin{2k}{k}}H_{q,k-1}(\{2^{r-1}\})\cdot [2k]_q! \left(
\frac{q^{k(k+1)}}{\prod_{j=1}^{2k}[j]_q}
-\frac{q^{(N+1)k}}{\prod_{j=N-k+1}^{N+k}[j]_q}
\right)\\
&=\sum_{k=1}^{N-1}\frac{(-1)^{k-1+r}q^{\binom{k+1}{2}+k}}{[k]_q^3\qbin{2k}{k}}
H_{q,k-1}(\{2^{r-1}\})
-\sum_{k=1}^{N-1}\frac{(-1)^{k-1+r}q^{-\binom{k}{2}+Nk+k}}{[k]_q^3\qbin{N}{k}\qbin{N+k}{k}}H_{q,k-1}(\{2^{r-1}\})\\
&=\sum_{k=1}^{N}\frac{(-1)^{k-1+r}q^{\binom{k+1}{2}+k}}{[k]_q^3\qbin{2k}{k}}
H_{q,k-1}(\{2^{r-1}\}) -\sum_{k=1}^{N}\frac{(-1)^{k-1+r}q^{-\binom{k}{2}+Nk+k}}{[k]_q^3\qbin{N}{k}\qbin{N+k}{k}}H_{q,k-1}(\{2^{r-1}\}).
\end{align*}
Finally, putting the LHS and the RHS back together, we have, for $r\geq 1$,
\begin{align*}
\sum_{n=1}^{N}&\frac{q^{rn}+q^{(r+1)n}}{[n]_q^{2r+1}}
=-\sum_{n=1}^{N}\frac{(-1)^{n+r}q^{\binom{n}{2}}(1+q^n)^2}{\qbin{2n}{n}}
\sum_{j=1}^{r}
\frac{(-1)^jq^{nj}}{[n]_q^{2j+1}}H_{q,n-1}(\{2^{r-j}\})\\
&+\sum_{k=1}^{N}\frac{(-1)^{k+r}q^{\binom{k+1}{2}+k}}{[k]_q^3\qbin{2k}{k}}
H_{q,k-1}(\{2^{r-1}\})-\sum_{k=1}^{N}\frac{(-1)^{k+r}q^{-\binom{k}{2}+Nk+k}}{[k]_q^3\qbin{N}{k}\qbin{N+k}{k}}H_{q,k-1}(\{2^{r-1}\})
\end{align*}
which is the desired finite form of \eqref{qMarkovAperyr} after replacing $r$ with $r+1$.
\hfill$\square$

%%%%%%%%%%%%%%%%%%%%%%%%%%%%%%%%%%%%%%%%%%%%%%%%%%%%%%%%%%%%%%%%%%%%%%%%

\end{document}